\documentclass[11pt]{article}
\usepackage{color}
\usepackage{amsmath, amssymb}
\usepackage{graphicx}
\usepackage{float}

\newcommand{\cqfd}{{\nobreak\hfil\penalty50\hskip2em\hbox{}
\nobreak\hfil $\square$\qquad\parfillskip=0pt\finalhyphendemerits=0\par\medskip}}
\newcommand{\R}{\mathbb{R}}

\newcommand{\eps}{ \varepsilon}

\newcommand{\pr}{{\bf \textit{Proof : }}}

\numberwithin{equation}{section}

\def\disp{\displaystyle}

\def\G{\gamma}

\def\l{\leq}
\def\g{\geq}
\def\q{\quad}

\def\i{\infty}

\def \B{\begin{equation}}
\def\bbb{\begin{array}{lll}}

\def\eee{\end{array}}
\def \E{\end{equation}}
\def\iii{\int_0^1}

\def\it{\indent}
\def\w{\widetilde}

\newtheorem{theorem}{Theorem}[section]
\newtheorem{definition}{Definition}[section]

\newtheorem{lemma}{Lemma}[section]

\newtheorem{remark}{Remark}[section]

\begin{document}
\title{ \bf
Large time behavior of entropy solutions to one-dimensional unipolar hydrodynamic model for semiconductor devices}

\author{\small
\textbf{Feimin Huang}\thanks{School of Mathematical Sciences, University of Chinese Academy of Sciences, Beijing 100049, China; Academy of Mathematics and System Sciences,  CAS, Beijing 100190, China. E-mail: fhuang@amt.ac.cn,}\quad
\textbf{Tianhong Li}\thanks{Hua Loo-Keng Key Laboratory of Mathematics,
	AMSS, CAS, Beijing 100190, China. E-mail: thli@math.ac.cn}\quad
\textbf{Huimin Yu}\thanks{Department of Mathematics,
Shandong Normal University, Jinan 250014, China,
E-mail: hmyu@amss.ac.cn }\quad
\textbf{Difan Yuan}\thanks{School of Mathematical Sciences, University of Chinese Academy of Sciences, Beijing 100049, China; Academy of Mathematics and System Sciences,
 CAS, Beijing 100190, China. E-mail: yuandf@amss.ac.cn}
}

\date{}

\bibliographystyle{plain}

\maketitle

\abstract{
  We are concerned with the global existence and large time behavior of entropy solutions to the one dimensional unipolar hydrodynamic model for semiconductors in the form of Euler-Poisson equations in a bounded interval.
  In this paper, we first prove the global existence of entropy solution by vanishing viscosity and compensated compactness framework. In particular, the solutions are uniformly bounded with respect to space and time variables by introducing modified Riemann invariants and the theory of invariant region. Based on the uniform estimates of density, we further show that the entropy solution converges to the corresponding unique stationary solution exponentially in time. No any smallness condition is assumed on the initial data and doping profile. Moreover, the novelty in this paper is about the unform bound with respect to time for the weak solutions of the isentropic Euler-Possion system.
\medskip

\noindent {\bf 2010 AMS Classification}: 35L60, 35L65, 35Q35.

\section{Introduction}
%

We consider a one-dimensional isentropic Euler-Poisson model for
semiconductor devices:
\begin{eqnarray}\label{semi}
\left\{ \begin{array}{llll}
\displaystyle  n_t+J_x=0,\\
\displaystyle J_t+\left({\frac{J^2}{n}}+p(n)\right)_x=nE-\frac{J}{\tau},\\
\displaystyle E_{x}=n-D(x),\\
\end{array}
\right.
\end{eqnarray}
in the region $\Pi_T:=(0,1)\times[0,T)$ for some fixed $T>0$, where $n\geq0, J,$ and $E$ denote the electron
density,
 electron current density, and the electric field, respectively.  $ E$ is generated by the Coulomb force of the particles and the pressure function is
 $p(n)=p_0n^\gamma$ with $p_0=\frac{\theta^2}{\gamma},\theta=\frac{\gamma-1}{2}$. Here $\gamma>
 1$ corresponds to the isentropic case. $\tau=\tau(n,J)>0$  is the momentum relaxation time.
Since we consider the large time behavior of solutions, rather than the relaxation limit, we assume $\tau=1$ for simplicity in this paper.
 The doping profile $D(x)\geq0$ stands for the density of fixed, positively charged
 background ions. For simplicity, we assume
 \begin{equation}\label{D(x)}
 \begin{aligned}
&D(x)\in C[0,1] ~{\rm{and }}~ 0<D_*\l D(x)\l D^*,
\end{aligned}
\end{equation}
where $D_*$ and $D^*$ are two positive constants, otherwise, we may assume that $D(x)$ is well approximated by a smooth function and treat it in a similar way as in \cite{G}.
The system $\eqref{semi}$ is supplemented by the following initial-boundary value conditions:
\begin{equation}\label{initial}
\begin{split}
&(n, J)|_{t=0}=(n_0(x), J_0(x)),0<x<1,\\
&J|_{x=0}=J|_{x=1}=0, t\geq0,\\
& E(0,t)=0, t\geq0,
\end{split}
\end{equation}
where $n_0(x)$ satisfies
\begin{equation}\label{n0}
\int_0^1 \big{(}n_0(x)-D(x)\big{)}dx=0.
\end{equation}
From $\eqref{semi}$ and $\eqref{initial},$ we have
 \begin{equation}\label{E}
 \begin{aligned}
E(x,t)=\int_0^x\left(n(\xi,t)-D(\xi)\right)d\xi.
\end{aligned}
\end{equation}
Then the system $\eqref{semi}$ can be reduced to the following one:
\begin{eqnarray}\label{semi2}
\left\{ \begin{array}{llll}
\displaystyle  n_t+J_x=0,\\
\displaystyle J_t+\left({\frac{J^2}{n}}+p(n)\right)_x=n\left(\int_0^x\left(n(\xi,t)-D(\xi)\right)d\xi\right)-{J},\\
\end{array}
\right.
\end{eqnarray}
which is a system of Euler equations in nonconservative form with inhomogeneous sources. Notice that the first source term, caused by the field $E$, is a nonlocal term involving some global properties of the solution of  $\eqref{semi2}$.
 For the detail derivation of the hydrodynamic model of semiconductors, one can see  \cite{Markowich} and \cite{Poupaud1}. \\
\it Now let's recall some known results for the problem
of \eqref{semi} in homogeneous case, i.e., one-dimensional Euler equation for compressible
fluid,
\begin{eqnarray}\label{homogeneous}
\left\{ \begin{array}{llll}
\displaystyle  n_t+J_x=0,\\
\displaystyle J_t+\left({\frac{J^2}{n}}+p(n)\right)_x=0,\\
\end{array}
\right.
\end{eqnarray}
which describes the motion of the ideal gas.
For the isentropic gas, the existence of BV solution was established by Nishida and
Smoller \cite{Smoller1} for some restricted classes of initial data.
It is well known that the system \eqref{homogeneous} is singular at $n=0$ (vacuum) and this poses new
difficulties to the mathematical analysis. When vacuum occurs, the
first global existence result of $L^\infty$ weak solution with large
initial data was firstly established by Diperna \cite{Diperna} for
$\gamma=1+\frac{2}{2n+1},\ n\geq2$, using the method of compensated
compactness. The existence problem was solved by Ding, Chen, and
Luo \cite{Ding1} for general $\gamma$ in the interval $(1,\frac{5}{3}]$, and further  extended to the case
$\gamma>{5\over 3}$ by Lions-Perthame-Tadmor \cite{Lions1} and
Lions-Perthame-Souganidis \cite{Lions2}, finally the isothermal case $\gamma=1$ by Huang-Wang \cite{HuangWang}. \\
\it Concerning the inhomogeneous case, i.e., the isentropic Euler equations with source term, Ding {\textit{et al}}.\cite{Ding1} established a general framework to investigate the global existence of entropy solution through the fractional step Lax-Friedrichs scheme and compensated compactness method. It should be noted that in the framework of \cite{Ding1}, the approximate solutions are only required to be uniformly bounded with respect to the space $x$, but not the time $t$. The $L^\infty$ norm of approximate solutions may increase with respect to time $t$. Later on, there have been extensive studies on the inhomogeneous case, see \cite{Chen2, Marcati, Marcati2} and the references therein. \\
\it We now turn to the Euler-Poisson system \eqref{semi}. A lot of
efforts have been made recently for this system on the whole space or bounded
domain. The
 mathematical study was initiated by Degond and Markowich \cite{Degond} in which  the existence of a unique smooth solution was proved for the steady-state of \eqref{semi} in the subsonic case, which is characterized by a smallness assumption on the current flowing through the device. The existence of a local smooth solution of the time-dependent problem was proved by using Lagrangian mass coordinates in \cite{Zhangbo}. As for the weak solutions,
Zhang \cite{Zhang} and Marcati and Natalini \cite{Marcati} investigated the
global existence of entropy solutions to the initial-boundary
value problem and the Cauchy problem for $\gamma>1$, respectively. Li \cite{Li} and Huang \textit{et al.} \cite{Huang1}
used a fractional Lax-Friedrichs scheme to prove the existence of
$L^\infty$ entropy solution of \eqref{semi} with $\gamma=1$ on a bounded
interval and the whole space. It is noted that the $L^\infty$ estimates of entropy solution in all of the above works depend on time $t$. More interesting model for one dimensional hydrodynamic model of semiconductor devices can be found in \cite{Gasser},\cite{Michele2013},\cite{Michele2017} and the references therein.\\
\it The large time behavior of entropy solution to the system (\ref{semi}) was first investigated in \cite{HPY} based on an assumption that  the uniform bound of density is independent of time $t$. Later in \cite{Yu}, the assumption of  \cite{HPY} was  relaxed to the case that the density can increase {\bf{slowly}} with time, i.e.,
\begin{eqnarray}\label{T}
\|n(x,t)\|_{L^\i}\l Ct^2,
\end{eqnarray}
where $\G>3$.

In this paper, we shall verify the assumption of \cite{HPY} that the density is uniformly bounded with respect to space $x$ and time $t$ in a bounded interval for $1<\G\le 3$. Different from usual ways as in \cite{Ding1,Marcati}, we construct the approximate solution by adding artificial viscosity perturbations to the system ({\ref{semi}}). Then we
introduce two modified Riemann invariants and derive a decoupled parabolic system
along characteristics. Due to the hyperbolicity of the system ({\ref{semi}}),
the integral of source terms along characteristics with respect to time $t$ can be translated to
the integral with respect to space $x$. Thus the uniform bounds of viscosity solutions independent of time $t$ can be derived through the theory of invariant region \cite{Smoller}.
Based on the uniform estimates of density,
we can further investigate  the global existence and large time behavior of entropy solutions. The idea can be applied to the other hyperbolic systems of conservation laws with source terms, see \cite{Huang2017} for the multi-dimensional compressible Euler equations with spherical symmetry.

 Before stating the main results, we define the entropy solution of the system \eqref{semi2} as follows.
\begin{definition}\label{def}
For every $T>0,$ a pair of bounded measurable functions $v(x,t)=(n(x,t),J(x,t))$ is called entropy solution of \eqref{semi2} with initial-boundary condition \eqref{initial} if
\begin{equation}\label{weak}
\begin{split}
&\int_0^T\int_0^1(n\psi_t+J\psi_x)dxdt+\int_{t=0}n_0\psi dx=0,\\
&\int_0^T\int_0^1\left(J\psi_t+\left(\frac{J^2}{n}+p(n)\right)\psi_x\right)dxdt\\
&+\int_0^T\int_0^1\left(n\left(\int_0^x\left(n(\xi,t)-D(\xi)\right)d\xi\right)-J\right)\psi dxdt+\int_{t=0}J_0\psi dx=0.
\end{split}
\end{equation}
for all $\psi\in C_0^\infty(\Pi_T)$,  $\Pi_T:=(0,1)\times[0,T)$
and the boundary condition is satisfied in the sense of divergence-measure field \cite{Chen1}, and for any weak and convex entropy pairs $(\eta,q)$, the entropy condition
\begin{equation}\label{entropy}
\begin{split}
\eta_t+q_x\l \eta_J(nE-J)
\end{split}
\end{equation}
is satisfied in the sense of distribution.
\end{definition}
The first result of the present paper is given below.
\begin{theorem}[\textbf{Existence and uniform estimates}]\label{main}
Let $1<\gamma\l 3.$ Assume that the initial data satisfy
\begin{equation}\label{ini}
0\leq n_0(x)\leq M_1,~~|J_0(x)|\leq n_0(x)M_1,~~\text{ a.e. },
\end{equation}
for some positive constant $M_1$.
Then there exists a global entropy solution of the initial-boundary value problem \eqref{semi}\eqref{initial} satisfying
\begin{equation}\label{solution}
0\leq n(x,t)\leq C,~|J(x, t)|\leq Cn(x, t),~|E(x, t)|\leq C,~~ \text{ a.e. },
\end{equation}
where $C$ only depends on $M_1$ and is independent of time $t$.
\end{theorem}
\begin{remark}
	Theorem \ref{main} is only valid for the bounded  interval in space. It is very interesting to show \eqref{solution} for the Cauchy problem.
	\end{remark}
	
	The second result of the present paper is given as follows.
\begin{theorem}[\textbf{Large time behavior}]\label{main2}
Denote $(n,J, E)$ is the global entropy solution of \eqref{semi}\eqref{initial} obtained in Theorem \ref{main} and  $(\w{N},\w{E})$ is the corresponding steady solution. Then it holds that
\begin{equation}\label{ini}
\iii \bigg{(}(n-\w{N})^2+(E-\w{E})^2+J^2\bigg{)}(x,t)dx\l C e^{-Ct}
\end{equation}
for some positive constant $C>0$.
\end{theorem}

Throughout this paper, $C_i, \w{C}_i (i=1,2\cdots)$ denote some specific positive constants, $C$ denotes the generic positive constant, which can be different from each other. The present paper is organized as follows. In Section \ref{formula}, we give some preliminaries of the system \eqref{semi} and approximate it by adding artificial viscosity.  In Section \ref{approximate}, we obtain the uniform bounds of viscosity solutions by introducing two modified Riemann invariants and the theory of invariant region.
Section \ref{convergence} is devoted to the proof of  Theorem \ref{main} and Theorem \ref{main2}.
\section{Preliminary and formulation}\label{formula}

We first introduce some basic facts for the system \eqref{homogeneous}.
The eigenvalues are
\begin{equation}\label{2.1}
\lambda_1=\frac{J}{n}-\theta n^\theta,\quad
\lambda_2=\frac{J}{n}+\theta n^\theta,\quad
\end{equation}
where $\theta=\frac{\gamma-1}{2},$ and the corresponding right eigenvectors are
\begin{equation}\label{2.2}
r_1=\left[\begin{array}{ll}
1\\ \lambda_1
\end{array}
\right],\quad
r_2=\left[\begin{array}{ll}
1\\ \lambda_2
\end{array}
\right].
\end{equation}
The Riemann invariants $(w, z)$ are given by
\begin{equation}\label{2.3}
w=\frac{J}{n}+n^\theta,\quad z=\frac{J}{n}-n^\theta,
\end{equation}
satisfying $\nabla w\cdot r_1=0$ and $\nabla z\cdot r_2=0,$ where $\nabla=(\partial_n, \partial_J)$ is the gradient with respect to $U=(n,J)$. A pair of functions $(\eta, q):\R\times\R^+\mapsto\R^2$ is called an entropy-entropy flux pair of system \eqref{homogeneous} if it satisfies
\begin{equation}\label{2.4}
\nabla q(U)=\nabla\eta(U)\nabla\left[\begin{array}{lll}
~~~~J\\ \frac{J^2}{n}+p(n)
\end{array}
\right].
\end{equation}
And $\eta(n, J)$ is called weak entropy if $$\eta\left|_{n=0, \frac{J}{n}\text{ fixed }}=0\right..$$
Moreover, an entropy $\eta(n,J)$ is said to be convex (strictly convex) if the Hessian matrix $\nabla^2\eta(n, J)$ is nonnegative (positive). For example, the mechanical energy
\begin{equation}\label{2.5}
\eta_e(n, J)=\frac{J^2}{2n}+\frac{p_0n^\gamma}{\gamma-1},~~
q_e(n, J)=\frac{J^3}{2n^2}+\frac{p_0n^{\gamma-1}J}{\gamma-1},
\end{equation}
is a strictly convex entropy pair.
As shown in \cite{Lions1} and \cite{Lions2},  any weak entropy for the system \eqref{homogeneous} is
\begin{equation}\label{2.6}
\begin{split}
\eta=n\int_{-1}^1g(\frac{J}{n}+n^\theta s)(1-s^2)^\lambda ds,~~
q=n\int_{-1}^1(\frac{J}{n}+n^\theta\theta s)g(\frac{J}{n}+n^\theta s)(1-s^2)^\lambda ds,
\end{split}
\end{equation}
with  $\lambda=\frac{3-\gamma}{2(\gamma-1)}$ and $g(\cdot)\in C^2(\R)$ is any function.
We approximate the system \eqref{semi2} by adding the following artificial viscosity, i.e.,
\begin{eqnarray}\label{semi-vis}
\left\{ \begin{array}{llll}
\displaystyle n_t+J_x=\eps n_{xx},\\
\displaystyle J_t+\left(\frac{J^2}{n}+p(n)\right)_x
=\eps J_{xx}+nE-{J}-2\eps n_x,
\end{array}
\right.
\end{eqnarray}
with initial-boundary data
\begin{equation}\label{ini-vis}
\begin{split}
&(n, J)|_{t=0}=(n_0^\eps(x), J_0^\eps(x))=(n_0(x)+\eps,J_0(x))\ast j^{\eps},0\leq x\leq1,\\
&{J}(0,t)={J}(1,t)=0,~n(0,t)=n_0^{\eps}(0),n(1,t)=n_0^{\eps}(1), t>0,
\end{split}
\end{equation}
where $n_0^{\eps}(x)\geq\eps>0,$ $j^{\eps}$ is the standard mollifier and the parameter $\eps>0$ is small.
We shall prove that the viscosity solutions of \eqref{semi-vis}-\eqref{ini-vis} are uniformly bounded with respect to time $t$  in the next section.

\section{Viscosity solutions and a priori estimates}\label{approximate}
\par
For any fixed $\eps>0$, we denote the solution of \eqref{semi-vis}-\eqref{ini-vis} by $(n^\eps, J^\eps)$.
The local existence of the approximate solution can be proved by the same argument of \cite{Diperna}. To extend the local solution to the global one, the key point is to
obtain the uniform bound of $n^\varepsilon, J^\eps$ and the lower bound of density $n^\eps$.

\subsection{Maximum principle}\label{upper}
\par
 We first give a modified maximum principle for a coupled quasilinear parabolic system, which is exactly a variant of the theory of invariant region \cite{Smoller}.
\begin{lemma}(Maximum principle)\label{maximum}
Let $(p, q)(x,t)$, $(x,t)\in  [0,1]\times[0,T]$ be any bounded classical solution of the following quasilinear parabolic system, i.e., $(p, q)(x,t)\in  C^{2,1}((0,1)\times[0,T])\cap C^0([0,1]\times[0,T])$ satisfy
\begin{eqnarray}\label{pq}
\left\{ \begin{split}
\displaystyle &p_t+\mu_1 p_x=
 p_{xx}+a_{11}p+a_{12}q+R_1,\\
\displaystyle &q_t+\mu_2 q_x=
 q_{xx}+a_{21}p+a_{22}q+R_2,
\end{split}
\right.
\end{eqnarray}
with initial-boundary data
$$p(x,0)\leq0, ~~q(x, 0)\geq0,$$$$p(0,t)\leq0,p(1,t)\leq0,q(0,t)\geq0,q(1,t)\geq0, $$
where $$\mu_{i}=\mu_i(x,t,p,q),a_{ij}=a_{ij}(x,t,p,q),$$ and the source terms $$R_i=R_i(x,t,p,q,p_{x},q_{x}),i,j=1,2,\forall(x,t)\in(0,1)\times[0,T].$$  Assume that $\mu_{i},a_{ij}$ are bounded with respect to $(x,t,p,q)\in[0,1]\times[0,T]\times K,$ where $K$ is an arbitrary compact subset in $\R^2,$  $a_{12},a_{21},R_{1},R_{2}$ are continuously differentiable with respect to $p,q$ and the following conditions hold:\\
$\bullet$ when $p=0$ and $q\geq0,$  $$a_{12}\leq0,~R_1=R_1(x,t,p,q,\zeta,\eta)\leq0;$$
$\bullet$ when $q=0$ and $p\leq0,$  $$a_{21}\leq0,~R_2=R_2(x,t,p,q,\zeta,\eta)\geq0.$$
Then for any $(x, t)\in[0,1]\times[0,T],$
$$p(x,t)\leq0, ~~q(x, t)\geq0.$$
\end{lemma}
\subsection{Uniform bound}\label{upper}
\par
Next we will derive the uniform bound of viscosity solution by the above maximum principle, i.e., Lemma \ref{maximum}.  We have
\begin{lemma}\label{uniform bound}
 If the viscosity solution $n^\varepsilon(x,t)\g 0$ for any $(x,t)\in (0,1)\times [0,T]$, then there exists a positive constant $C>0$,  independent of $T$ and $\eps$, such that
\begin{equation}\label{bound}
0\le n^\varepsilon(x, t)\leq C,~~ |J^\eps(x, t)|\leq Cn^\eps(x, t), ~(x,t)\in (0,1)\times [0,T].
\end{equation}
\end{lemma}
\pr
\footnote{For simplicity of notation, the superscript of $(n^\eps, J^\eps)$ will be omitted as $(n, J)$.}
By the formulas of Riemann invariants \eqref{2.3}, we can decouple the viscous perturbation equations \eqref{semi-vis} as
\begin{eqnarray}\label{wz}
\left\{ \begin{split}
\displaystyle &w_t+\lambda_2 w_x=\eps w_{xx}+2\eps(w_x-1)\frac{n_x}{n}-\eps\theta(\theta+1)n^{\theta-2}n_x^2+\frac{1}{n}\left(nE-{J}\right),\\
\displaystyle &z_t+\lambda_1 z_x=\eps z_{xx}+2\eps(z_x-1)\frac{n_x}{n}+\eps\theta(\theta+1)n^{\theta-2}n_x^2+\frac{1}{n}\left(nE-{J}\right).
\end{split}
\right.
\end{eqnarray}
Motivated by \cite{Tsuge2} in which a modified numerical scheme is used, set the control functions $(\varphi,\psi)$ as
\begin{equation*}
\begin{split}
&\varphi=M+x,\\
&\psi=M-x,
\end{split}
\end{equation*}
where $M$ is a constant to be determined.
A direct calculation tells us
\begin{equation*}
\begin{split}
&\varphi_t=0,~\varphi_x=1, ~\varphi_{xx}=0;\\
&\psi_t=0,~\psi_x=-1, ~\psi_{xx}=0.
\end{split}
\end{equation*}
Define the modified Riemann invariants $(\bar{w},\bar{z})$ as
\begin{equation}\label{r}
\bar{w}=w-\varphi, ~~\bar{z}=z+\psi.
\end{equation}
The advantage of using the modified Riemann invariants is that we can use the derivative of control function $\phi_x,\psi_x$ and eigenvalues $\lambda_1$ and $\lambda_2$ in the left hand side of \eqref{wz} so that the coefficients $a_{12}$, $a_{21}$ and the terms $R_1$ and $R_2$ are designed to have desired signs, see \eqref{sign} below.  Inserting the above formulas into \eqref{wz} yields the decoupled equations for $\bar{w}$ and $\bar{z}:$
\begin{eqnarray}\label{phipsi1}
\left\{ \begin{split}
\bar{w}_t+\lambda_2\bar{w}_x&=\eps\bar{w}_{xx}
+2\eps\bar{w}_x\frac{n_x}{n}-\eps\theta(\theta+1)n^{\theta-2}n_x^2+\frac{1}{n}\left(nE-{J}\right)-\lambda_2\varphi_x,\\
\bar{z}_t+\lambda_1\bar{z}_x&=\eps\bar{z}_{xx}
+2\eps\bar{z}_x\frac{n_x}{n}+\eps\theta(\theta+1)n^{\theta-2}n_x^2+\frac{1}{n}\left(nE-{J}\right)+\lambda_1\psi_x.
\end{split}
\right.
\end{eqnarray}
Note that
\begin{equation}
\begin{split}
&\lambda_1=\frac{w+z}{2}-\theta\frac{w-z}{2},\\
&\lambda_2=\frac{w+z}{2}+\theta\frac{w-z}{2},
\end{split}
\end{equation}
and then the system \eqref{phipsi1} becomes
\begin{eqnarray}
\displaystyle\left\{ \begin{split} &\bar{w}_t+(\lambda_2-2\eps\frac{n_x}{n})\bar{w}_x
=\eps\bar{w}_{xx}+a_{11}\bar{w}
+a_{12}\bar{z}+R_1,\\
&\bar{z}_t+(\lambda_1-2\eps\frac{n_x}{n})\bar{z}_x
=\eps\bar{z}_{xx}+a_{21}\bar{w}
+a_{22}\bar{z}+R_2.
\end{split}
\right.
\end{eqnarray}
with
\begin{equation}\label{coefficient}
\begin{split}
&a_{11}=-1-\frac{\theta}{2},~~a_{12}=-1+\frac{\theta}{2},~~a_{21}=-1+\frac{\theta}{2},~~a_{22}=-1-\frac{\theta}{2},\\
&R_1=-\eps\theta(\theta+1)n^{\theta-2}n_x^2+E-\left(1+\frac{\theta}{2}\right)\varphi+\left(1-\frac{\theta}{2}\right)\psi,\\
&R_2=\eps\theta(\theta+1)n^{\theta-2}n_x^2+E-\left(1-\frac{\theta}{2}\right)\varphi+\left(1+\frac{\theta}{2}\right)\psi.
\end{split}
\end{equation}
We shall use Lemma \ref{maximum} to show $\bar{w}\le 0$ and $\bar{z}\ge 0$. \\
\it Define a smooth function $\phi_\delta, $ such that for small $\delta>0$:
\begin{eqnarray}\label{phidelta}
\displaystyle\phi_\delta=\left\{ \begin{split} &x^2,~~~~~~ 0<x<\delta,\\
&1, ~~~~~2\delta<x<1-2\delta,\\
&(1-x)^2, ~~~1-\delta<x<1,\\
\end{split}
\right.
\end{eqnarray}
and for all $x\in [0,1]$,
\begin{equation}\label{phidelta2}
|\phi'_\delta|\leq\frac{C}{\delta},~|\phi''_\delta|\leq \frac{C}{\delta^2}.
\end{equation}
Multiplying the first equation in \eqref{semi-vis} by $\phi_\delta$ and  using the boundary condition \eqref{ini-vis}, we get
\begin{equation}\label{conservation}
\left(\int_0^1n\phi_\delta dx\right)_t-\int_0^1n u\phi'_\delta dx=\eps\int_0^1n\phi_\delta''dx.
\end{equation}
It is noted that
\begin{equation}\label{phidelta3}
\int_0^1nu\phi'_\delta dx=\int_0^\delta 2nu xdx-\int_{1-\delta}^12nu(1-x)dx+\int_\delta^{2\delta}nu\phi'_\delta dx+\int_{1-2\delta}^{1-\delta} nu\phi'_\delta dx.
\end{equation}
The first two terms in the right hand side of the above equality tend to zero when $\delta$ is small. The last two terms also tend to zero  by using the boundary condition \eqref{ini-vis} and \eqref{phidelta2} as $\delta$ is small enough.
Suppose that \begin{equation}\label{assu} n\le M_1=(2M)^{\frac1{\theta}}.
\end{equation} Then a simple calculation gives that
\begin{equation}\label{rr}
\left|\eps\int_0^1n\phi_\delta'' dx\right|\leq C\eps\int_0^{1}\frac{n}{\delta^2} dx
\leq CM_1\frac{\eps}{\delta^2}.
\end{equation}
Choosing $\delta=\eps^{\frac14}$. When $\eps$ is small,
we have $\left|\eps\int_0^1n\phi_\delta'' dx\right|\leq \frac12$ which gives that
 \begin{equation}\label{assu1}\left|\int_0^1n(x,t)\phi_{\delta} dx\right|\leq\left|\int_0^1n_0(x)\phi_{\delta} dx\right|+1.\end{equation}
Then,
\begin{equation}\label{n}
\begin{split}
\left|\int_0^1n dx\right|\leq& \left|\int_0^{1}n(\phi_{\delta}-1) dx\right|+\left|\int_0^{1}n\phi_{\delta} dx\right|\\
\leq&\left|\int_0^{2\delta}n(\phi_{\delta}-1) dx+\int_{1-2\delta}^{1}n(\phi_{\delta}-1) dx\right|+\left|\int_0^1n_0dx\right|+1\\
\leq& 8\delta M_1+\left|\int_0^1n_0dx\right|+1.
\end{split}
\end{equation}
 From \eqref{D(x)} and \eqref{assu1}, there exists a positive constant $M_0$ independent of $M_1$  such that,
$$|E(x,t)|\l \int_0^1 n(x,t)dx+\int_0^1 D(x)ds\leq M_0.$$
Hence, we can get from \eqref{coefficient} that
\begin{equation}\label{r1}
R_1\leq M_0-\left(1+\frac{\theta}{2}\right)(M+x)+\left(1-\frac{\theta}{2}\right)(M-x)\leq M_0-\theta M,
\end{equation}
and
\begin{equation}
\label{r2}
R_2\geq-M_0-\left(1-\frac{\theta}{2}\right)(M+x)+\left(1+\frac{\theta}{2}\right)(M-x)
\geq-M_0+\theta M-2.
\end{equation}
Noting that $0<\theta\leq 1,$ and choosing $M$ large enough, one has
\begin{equation}\label{sign}
a_{12}\leq0,a_{21}\leq0,R_1\leq0,R_2\geq0.
\end{equation}
Note also that on the boundary conditions,
\begin{equation}\bar{w}(0,t)\leq0,\bar{z}(0,t)\geq0,\bar{w}(1,t)\leq0,\bar{z}(1,t)\geq0,\end{equation}
and the initial values satisfy
 \begin{equation}\begin{split}\bar{w}(x,0)=w(x, 0)-\phi(x, 0)=\frac{J_0}{n_0}+n_0^\theta-M-x\leq0,\\\bar{z}(x,0)=z(x, 0)+\psi(x, 0)=\frac{J_0}{n_0}-n_0^\theta+M-x\geq0,
 \end{split}
 \end{equation}
 thanks to  Lemma \ref{maximum}, we have
$$\bar{w}(x, t)\le0,~~ \bar{z}(x, t)\ge0, ~\forall (x,t)\in[0,1]\times[0,T],$$
which  implies that
\begin{equation}\label{estimate1}
\begin{split}
&w(x, t)\leq\phi(x, t)\leq M+x\leq M+1,\\
&z(x,t)\geq-\psi(x,t)\geq-M+x\geq-M.
\end{split}
\end{equation}
By \eqref{estimate1}, we have
\begin{equation}
n\le (\frac{w-z}{2})^{\frac1{\theta}}\le (\frac32M)^{\frac1{\theta}}
\end{equation}
which verifies the assumption \eqref{assu}.  Therefore, Lemma \ref{upper} is completed.

From \eqref{bound}, the velocity $u=\frac{J}{n}$ is uniformly bounded, i.e., $|u|\le C$. Then following the same way of \cite{Huang2017},
we can obtain
\begin{equation}\label{density}
n(x,t)\geq\delta(t,\eps)>0,
\end{equation}
where $\delta(t,\eps)$ depends on $n_0,\eps, T$. Based on the local existence of solution, the uniform upper estimates (Lemma \ref{uniform bound}) and the lower bound estimate of density (\ref{density}), we can conclude the following theorem.
\begin{theorem}\label{ex}
For any time $T>0$, there exists a unique global classical solution $(n^\eps, J^\eps)(x,t), (x.t)\in [0,1]\times [0,T)$ to the initial-boundary value problem \eqref{semi-vis}-\eqref{ini-vis} satisfying
\begin{equation}\label{boundvis}
0<\delta(t, \eps)\leq n^\eps(x, t)\leq C, ~~|J^\eps(x, t)|\leq Cn^{\eps}(x, t),
\end{equation}
where $C$ is independent of $\eps$ and $T$.
\end{theorem}

 By Theorem \ref{ex} and the compactness framework established in
 \cite{Ding, Diperna, Lions2, Marcati, Zhang},  we can prove that there exists
 a subsequence of $(n^\eps,J^\eps)$ (still denoted by $(n^\eps,J^\eps)$) such  that
 \begin{equation}
 \label{4.6}
 (n^\eps, J^\eps)\to(n, J) ~~~
 \text{ in } L^p_{loc}(\Pi_T), ~~p\geq1. \end{equation}
 As in \cite{Chen2, Tsuge2},  we can prove that $(n, J)$ is an entropy  solution to the initial-boundary value problem \eqref{semi2}. The proof of Theorem \ref{main} is completed.
\section{Large time behavior of weak solutions}\label{convergence}

In this section, we shall show the entropy solution obtained in Theorem \ref{main} exponentially converges in time to the corresponding stationary solution, i.e. the smooth solution of
\begin{eqnarray}\label{steady}
\left\{ \begin{split}
&p(\w{N})_x=\w{N}\w{E},\\
&\w{E}_x=\w{N}-D(x),\\
\end{split}
\right.
\end{eqnarray}
with the boundary condition
\begin{equation}\label{B}
\w{E}(0)=\w{E}(1)=0.
\end{equation}
From ($\ref{steady}$) and ($\ref{B}$), we know that
$$
\iii \Big{(}\w{N}(x)-D(x)\Big{)} dx=0.
$$
The work of \cite{HY} told us the following result for the stationary solution:
\begin{lemma}\label{steady solution}{\rm{(Hsiao-Yang, 2001)}}
Under the assumption \eqref{D(x)} on $D(x)$, there exists a unique solution $(\w{N}, \w{E})$ to problem \eqref{steady} and \eqref{B} satisfying
\begin{equation}\label{bs}
0<D_*\l \w{N}(x)\l D^*,~ |\w{N}'(x)|\l C, ~|\w{N}''(x)|\l C, ~x\in [0,1],
\end{equation}
\end{lemma}
where $C$ only depends on $\G$, $D^*$ and $D_*$.
\cqfd

\it Based on the uniform estimates of viscosity solutions (Theorem \ref{main}), we shall prove that the entropy solution $(n, J, E)$ obtained in Theorem \ref{main} is strongly convergent to the corresponding stationary solution $(\w{N}, \w{E})$ in $L^2(0,1)$ norm with exponential decay rate. To this end, define new variable
\begin{equation}
y(x,t)=-\int_0^x\Big{(}n (s,t)-\w{N}(s)\Big{)}ds=-(E-\tilde{E}),\q (x,t)\in [0,1]\times [0,+\i).
\end{equation}
It is obvious that $y$ is absolutely continuous in $x$.
Then
\B\label{eq:44}
y_x=-(n -\tilde{N }),\q y_t=J.\E
From \eqref{n0}, $y(0)=y(1)=0$. The systems (\ref{semi}) and  (\ref{steady}) infer that
$y$ satisfies
\begin{equation}\label{eq:45}
y_{tt}+(\frac{J^2}{n})_x+(p(n)-p(\w{N}))_x+y_t=-\w{N}y-\w{E}y_x+yy_x.
\E Multiplying $y$ with $(\ref{eq:45})$ and integrating over $(0,
1)$,
we get
\begin{equation}\label{basic}
\begin{array}{ll} &\disp {d\over {dt}}\iii
(yy_t+\frac12y^2)dx+\iii
\bigg{(}(p(n)-p(\w{N}))(n-\w{N})+(\w{N}-\frac{\w{E}_x}{2})y^2\bigg{)}\
dx\\\\
&\q \disp\l \iii y_t^2dx +\iii \frac{y_t^2}{n}y_x\ dx=\iii \frac{\w{N}}{n}y_t^2dx.\end{array}
\end{equation}
Thanks to (\ref{bs}), the Lemma 3.1 of \cite{HP} says that there exist positive constants $\w{C}_1$ and $\w{C}_2$ such that
\begin{equation}\label{29}
\w{C}_2(n-\w{N})^2\ge (p(n)-p(\w{N}))(n-\w{N})\g \w{C}_1(n-\w{N})^2= \w{C}_1y_x^2.
\end{equation}
Then \eqref{basic} implies that
\begin{equation}\label{eq:47}\bbb
 \disp{d\over {dt}}\iii
 (yy_t+\frac12y^2)dx+\w{C}_1\disp\iii y_x^2dx\disp+D_*\iii y^2dx\l\iii {\frac{\tilde{N }}n }y_t^2dx.\eee
\end{equation}

On the other hand, denote the relative entropy-entropy flux  by
\begin{equation}\label{31}
\begin{array}{ll}
&\disp{\eta}_*=\eta_e-{{p_0\w{N}^\gamma}\over {\G-1}}-{{p_0\gamma}\over {\gamma-1}}\w{N}^{\gamma-1}(n-\w{N}),\\\\
&\disp {q}_*=q_e -{{p_0\G}\over {\G-1}}\w{N}^{\gamma-1}J.
\end{array}
\end{equation}
The entropy inequality \eqref{entropy} indicates that the following entropy
inequality holds in the sense of distribution:
\begin{equation}\label{eq:49}\bbb
\eta^*_t+q_{x}^*&=&\disp\eta_{et}+q_{ex}-\frac{p_0\G}{\G-1}\tilde{N }^{\G-1}(n -\tilde{N })_t-\frac{p_0\G}{\G-1}(\tilde{N }^{\G-1}J)_x\\
&\l&
\disp JE-\frac{J^2}n -\frac{p_0\G}{\G-1}\tilde{N }^{\G-1}(n -\tilde{N })_t-\frac{p_0\G}{\G-1}(\tilde{N }^{\G-1}J)_x\\
&=&\disp JE-\frac{J^2}n-p_0\G\tilde{N}^{\G-2}y_t\tilde{N}_x.\eee
\end{equation}
 Noticing that
\begin{equation}\bbb
JE&=&\disp J\tilde{E}+J(E-\tilde{E})=y_t\tilde{E}-yy_t\\\\
&=&\disp {{p(\tilde{N})_x}\over \tilde{N}}y_t-yy_t
=\disp p_0\G\tilde{N}^{\G-2}y_t\tilde{N}_x-yy_t \eee
\end{equation}
and using the theory of divergence-measure fields \cite{Chen1}, we
have
\begin{equation}\label{eq:50}
{d\over{dt}} \iii (\eta^*+{1\over 2}y^2)dx+\iii{{y_t^2}\over n }dx\l
0.
\end{equation}
\it Let $\Lambda$ big enough such that $\Lambda> D^*+\|n\|_{L^\i}+1$.
Multiply (\ref{eq:50}) by  $\Lambda$,
and add the result to (\ref{eq:47}), we have
\begin{equation}\label{ineq}
{d\over {dt}}\iii (\Lambda \eta_*+{{\Lambda y^2}\over 2}+yy_t+{{y^2}\over 2})dx+\w{C}_2\iii (y^2+y_x^2+{y_t^2\over n})dx\l 0.
\end{equation}
Since $\|n(x,t)\|_{L^\i}\l C$ and
$$\eta_*\sim  y_x^2+{y_t^2\over n},$$
it is straightforward to see that
$$\Lambda \eta_*+{{\Lambda y^2}\over 2}+yy_t+{{y^2}\over 2}\sim  y^2+y_x^2+{y_t^2\over n}.$$
Then from \eqref{ineq}, the Gronwall inequality implies  Theorem \ref{main2}.
\section*{Acknowledgments}
Feimin Huang is partially supported by National Center for Mathematics and Interdisciplinary Sciences, AMSS, CAS and NSFC Grant No. 11371349 and 11688101. Huimin Yu is supported in part by (NSFC Grant No.11671237) and
Natural Science Foundation of Shandong Province (Grant No. ZR2015AM001). The research of Difan Yuan is supported by China Scholarship Council No. 201704910503.
%

\end{document}